\documentclass{amsart}

\newtheorem{theorem}{Theorem}[section]
\newtheorem{lemma}[theorem]{Lemma}

\theoremstyle{definition}

\theoremstyle{remark}

\numberwithin{equation}{section}

\newcommand{\set}[1]{\left\{#1\right\}}
\begin{document}

\title{Level sets of multiple ergodic
averages}

\author{Ai-Hua Fan}
\address{LAMFA, UMR 6140 CNRS, Universit\'e de Picardie,
33 rue Saint Leu, 80039 Amiens, France}
%\curraddr{}
\email{ai-hua.fan@u-picardie.fr}

\author{Lingmin Liao}
\address{LAMA,  CNRS UMR 8050, Universit\'e Paris-Est, 61 Avenue du G\'en\'eral de Gaulle, 94010, Cr\'eteil Cedex, France}
\email{lingmin.liao@u-pec.fr}

\author{Ji-Hua Ma}
\address{Department of Mathematics, Wuhan University, 430072 Wuhan, China}
\email{jhma@whu.edu.cn}

%    General info

%\date{January 1, 2001 and, in revised form, June 22, 2001.}



\begin{abstract}
We propose to study  multiple ergodic averages from multifractal
analysis point of view. In some special cases in the symbolic
dynamics, Hausdorff dimensions of the level sets of multiple ergodic
average limit are determined by using Riesz products.
\end{abstract}

\maketitle

\section{Introduction}
Let $(X,T)$ be a topological dynamical system and let $\ell\ge 2$ be
a positive integer. We consider the following multiple ergodic
averages
      \begin{equation}\label{MEA}
\frac{1}{n}\sum_{k=1}^{n}f_{1}(T^{k}x)f_{2}(T^{2k}x)\cdots
f_{\ell}(T^{\ell k}x),
      \end{equation}
where $f_1,\cdots,f_\ell$ are $\ell$ given continuous functions.
Such multiple ergodic averages were introduced and studied by
Furstenberg \cite{Furs77} in his ergodic theoretic proof of
Szemer\'{e}di's theorem on arithmetic progressions. Since then these
averages have received extensive studies in various contexts. For
example, the $L^2$-normal convergence of (\ref{MEA}) is proved by
Host and Kra \cite{HK} with respect to a given invariant measure,
and the almost sure convergence is  proved earlier by Bourgain
\cite{Bou} in the case of $\ell=2$. In this note we propose to study
these multiple ergodic averages from multifractal analysis point of
view.
%We should say that it is far from
%evident to adapt the classical methods and techniques in the
%multiple case, i.e. $\ell\ge 2$.

Multifractal analysis of ergodic averages concerns the Hausdorff dimension of the level sets of the ergodic average limit. It reflects the complex  
behavior of the underlying chaotic dynamical system. There was a wide study in the case of simple  
ergodic averages  ($\ell =1$) in the last decades (\cite{FF,FFW,FLP,MW, Oli, Pes, Sch}). Our  
first investigation shows that the multifractal analysis  of multiple  
ergodic averages ($\ell\ge 2$) is much more difficult. This note aims at a special case where $X$ is the symbolic space
$\mathbb{D}=\set{+1,-1}^{\mathbb{N}}$ ($\mathbb{N}$ denoting the set
of positive integers) and the dynamics is defined by the shift
transformation $T:(x_1,x_2,\cdots)\mapsto(x_2,x_3\cdots)$. The
metric on $\mathbb{D}$ is chosen to be
  $$\rho(x,y)=2^{-\min\{k\geq 1:x_{k}\not=y_{k}\}} \ \textrm{ for } x,y\in \mathbb{D}.$$
The Hausdorff dimension of a set $A$ will be denoted by $\dim_H A$.
See \cite{Fal} for notions of dimensions of a set and \cite{Fan1994}
for notions of dimensions of a measure.
Let $\ell\geq 1$. We shall examine the averages (\ref{MEA}) with
the functions
  \begin{equation}\label{choice}
  f_1(x)=f_2(x)=\cdots=f_{\ell}(x)=x_1 \ \textrm{ for } x\in \mathbb{D}.
  \end{equation}
Then  for $\theta\in[-1,1]$, we consider the level set
            \begin{equation*}
B_{\theta}:=\set{x\in\mathbb{D}:
    \lim_{n\to\infty}\frac{1}{n}\sum_{k=1}^{n} x_kx_{2k}\cdots x_{\ell k}=\theta}.
       \end{equation*}

We prove the following result. \medskip
   \begin{theorem}\label{Main} For any $\theta\in[-1,1]$, we have
$$\dim_H(B_{\theta})=1-\frac{1}{\ell}+\frac{1}{\ell}H\Bigl(\frac{1+\theta}{2}\Bigr),$$
where $H(t)=-t\log_{2} t-(1-t)\log_{2}(1-t)$ is the entropy
function.
    \end{theorem}
    \medskip

This result was known to Besicovitch and Eggleston when $\ell=1$.
 Remark that the Hausdorff dimension of $B_{\theta}$
is strictly positive for any $\theta\in[0,1]$ when $\ell\geq 2$.
Actually,
$$ \dim_H B_\theta \ge 1 -1/\ell
>0\quad \mbox{ if} \quad \ell \ge 2. $$

The proof of the theorem is based on the fact that $\mathbb{D}$ has
a group structure and the functions $x\mapsto x_kx_{2k}\cdots
x_{\ell k}$ are group characters and even they  constitute a
dissociated set of characters in the sense of Hewitt-Zuckermann
\cite{HZ}. As we shall show, the set $B_\theta$ supports  a Riesz product, a nice
measure which has the same Hausdorff dimension as
that of $B_\theta$. 
The idea of using Riesz product is inspired by
\cite{Fan}  where oriented walks were studied. Although the Riesz product works  
perfectly for the above
case concerned by Theorem \ref{Main}, it has its limit for the general case.

We point out that the situation seems very different when the functions in (\ref{choice}) are replaced by other
functions. For
example, when $f_i$ are chosen as $(x_1+1)/2$ which takes $0$ and
$1$ as values. The obtained set can be identified with
     $$A_{\theta}:=\set{x\in \{0,1\}^{\mathbb{N}} :
    \lim_{n\to\infty}\frac{1}{n}\sum_{k=1}^{n} x_kx_{2k}\cdots x_{\ell k}=\theta}.$$
The set $A_\theta$  is  similar to $B_{\theta}$, but the
determination of its dimension is more difficult.

Actually, we are motivated by the study of $A_{\theta}$. The Riesz
product method is not adapted to it. Then we propose to looking at the following set
%which admits the same
%dimension as $A_0$ :
$$X_{0}:=\set{x\in \{0,1\}^{\mathbb{N}}: \ x_nx_{2n}=0, \quad \text{for all} \ n},$$
%We show that the Hausdorff dimensions of $A_0$ and $X_0$ are the same.
%\begin{thm}\label{Main-HD-same} $$\dim_H(A_0)=\dim_H(X_0).$$
%    \end{thm}
which is a subset of $A_{\theta}$ with $\ell=2$ and
$\theta=0$. We obtain the
box dimension (denoted by $\dim_B$) for $X_0$ by a combinatoric method.
\medskip
\begin{theorem}\label{Main-box} 
Let $\set{a_n}$ be the Fibonacci sequence
defined by
  $$a_0=1, \ a_1=2,\  \ a_{n}=a_{n-1}+a_{n-2} \ (n\geq 2).$$
We have
$$\dim_B(X_0)=\frac{1}{2\log 2}\sum_{n=1}^{\infty}\frac{\log a_n}{2^n}=0.8242936\cdots.$$
 \end{theorem}

The problem of determining the Hausdorff dimension of $X_0$ is now
solved by Kenyon, Peres and Solomyak \cite{KPS}, where a class of
sets similar to $X_0$ is studied. The result in \cite{KPS} together
with Theorem \ref{Main-box} shows that $\dim_HX_0 < \dim_B X_0$.

%One might guess this should be the Hausdorff dimension formula for both of $A_0$ and $X_0$. But this is not the case by the recent work of
%Kenyon, Peres and Solomyak \cite{KPS}, where the Hausdorff dimension of the set $X_0$ is obtained.

\medskip
\section{Riesz products}
  Let us consider $\mathbb{D}$ as an infinite product group of the multiplicative group
$\{+1,-1\}$. The dual group of $\mathbb{D}$ consists of the Walsh
functions $\set{w_n(x)}_{n=0}^{\infty}$ defined as follows. Define
$w_0=1$. For each $n\geq 1$, let
$$n=2^{n_1-1}+2^{n_2-1}+\cdots+2^{n_s-1},\quad  1\leq n_1<n_2<\cdots<n_s,$$
be the unique expansion of the integer $n$ in base $2$. Then we
define
    \begin{equation*}\label{Walsh}
w_n(x)=x_{n_1}x_{n_2}\cdots x_{n_s}.
    \end{equation*}
An important subset of Walsh functions is the set of the Rademacher functions
$\{r_n(x)\}_{n=1}^{\infty}$ defined by $r_n(x)=x_n$. The Rademacher
functions are mutually independent with expectation zero with respect to
the Haar measure. The following immediate consequence of the
independence will be frequently used in the sequel.
\medskip
      \begin{lemma}\label{Orth}  Let $f$ and $g$ be two Haar integrable functions on $\mathbb{D}$. Suppose that $f$ depends only on
the first $n$ coordinates of $x$ and $g$ is independent of the first
$n$ coordinates. Then 
             $$\int f(x)g(x)dx=\int f(x)dx\int g(x)dx$$
             where $dx$ stands for the Haar measure on $\mathbb{D}$.
%In particular, $\int f(x)g(x)dx=0$ if $g(x)=x_{n+k}$($k\geq 1$).
      \end{lemma}

\medskip
The $n$-th Fourier coefficient of an integrable function $f$ is defined by
\[\hat{f}(n)= \int f(x) w_n(x) dx.\]
   In the
follows, we shall denote
$$\xi_k(x)=x_kx_{2k}\cdots x_{\ell k} \ \ \ \textrm{ for all } k\geq 1.$$
   Consider the  product
     \begin{equation*}\label{Riesz}
dP_{\theta}(x)=\prod_{k=1}^{\infty}\bigl(1+\theta \xi_k(x)\bigr)dx.
     \end{equation*}
 The
following lemma shows that the above product defines a probability
measure on $\mathbb{D}$, which  will be called Riesz product.
\medskip

  \begin{lemma}\label{expectation} The partial products of the above infinite product
  converge in the weak-$*$ topology to a probability measure
  $P_{\theta}$. Furthermore,
   for any function $f$ depending only on the first $n$ coordinates of $x$, we have
             \begin{equation}\label{experence}
\mathbb{E}_{\theta}[f]=\int f(x)\prod_{k=1}^{\lfloor
n/{\ell}\rfloor}\bigl(1+\theta \xi_k(x)\bigr)dx,
             \end{equation}
where $\mathbb{E}_{\theta}[\;\cdot\;]$ stands for the expectation
with respect to $P_{\theta}$ and `` $\lfloor\;\cdot\; \rfloor$" is
the integer part function.
   \end{lemma}
\medskip
{\bf Proof}\quad For $N\geq 1$, let $$
P_N(x)=\prod_{k=1}^{N}\bigl(1+\theta
\xi_k(x)\bigr).$$
 Then
     $$  P_{N+1}(x)-P_N(x)=\theta P_N(x)\xi_{N+1}(x).$$
Observe that for the fixed Walsh function $w_n(x)=x_{n_1}x_{n_2}\cdots
x_{n_s}$, by Lemma \ref{Orth}, one has $$\int
P_N(x)\xi_{N+1}(x)w_n(x) d x=0$$ whenever $(N+1)\ell>n_s$. It
follows that %the $n$-th Fourier coefficient of $P_N$
$\hat{P}_N(n)=\hat{P}_{N+1}(n)$ for large $N$, so the limit
  $$\lim_{N\to\infty}\int P_N(x)w_n d x $$
exists. That is to say, the measures $P_N(x)dx$ converge weakly to a
limit measure $P_{\theta}$.

The formula (\ref{experence}) follows directly from Lemma \ref{Orth}
and the definition of the Riesz product $P_{\theta}$ as a weak limit. $\Box$
\medskip

The functions $\xi_n$ are not $P_\theta$-independent, but they are
orthogonal. %Even functions of $\xi_n$ are orthogonal. 
Therefore, we
can get the following law of large numbers.
  \begin{lemma}\label{Taylor} Suppose that $g$ is a function on the interval $[-1,1]$
  such that
  \[g(t)= \sum_{n=0}^{\infty} g_n t^n \quad \textrm{ with } \quad \sum_{n=1}^{\infty} |g_n|<\infty.\]
Then for $P_{\theta}$-almost all $x$,
$$\lim_{n\to\infty}\frac{1}{n}\sum_{k=1}^{n}g(\xi_k(x))=\mathbb{E}_{\theta}[g(\xi_1)].$$
  \end{lemma}
{\bf Proof}\quad Notice that $\xi_k^{2n}(x)=1 $ and
$\xi_k^{2n-1}(x)=\xi_k(x)$ for any integer $n\ge 1$. Then we get
   $$g(\xi_k)=\sum_{n=0}^{\infty} g_{2n} + \xi_k
\sum_{n=1}^{\infty} g_{2n-1}.$$ By the formula (\ref{experence}), we
have
\begin{equation*}\label{experence-variation}
\mathbb{E}_{\theta}(\xi_k)=\theta,\quad
\mathbb{E}_{\theta}(\xi_j\xi_k)=\theta^2,\quad(j\not=k).
\end{equation*} It follows that
$$
  \mathbb{E}_{\theta}[g(\xi_k)]=\sum_{n=0}^{\infty} g_{2n} + \theta
\sum_{n=1}^{\infty}
  g_{2n-1}, \qquad
\textrm{Cov}_{\theta} [g(\xi_j),g(\xi_{k})]=0 \quad(j\not=k).
$$
Therefore, the system
   $g(\xi_k)-\mathbb{E}_{\theta}[g(\xi_k)]$ ($k=1, 2, \cdots)$
is orthogonal in  $L^2(P_{\theta})$. By the Menchoff Theorem
(\cite{Zy}), the series
\[ \sum_{k=0}^{\infty} \frac{1}{k}\Big(g(\xi_k)-\mathbb{E}_{\theta}[g(\xi_k)]\Big)\]
 converges $P_{\theta}$-almost surely. Now the desired result follows from Kronecker's theorem.
 $\Box$

\medskip
\section{Proof of Theorem \ref{Main}}

Applying Lemma \ref{Taylor}  to $g(t)=t$, we get that for
$P_{\theta}$-almost all $x$,
\begin{equation*}
  \lim_{m\to \infty}\frac{1}{m}\sum_{k=1}^{m} \xi_k(x) =
  \mathbb{E}(\xi_1)=\theta.
\end{equation*}
This means that the Riesz product $P_{\theta}$ is supported by the
set $B_{\theta}$. Now we are going to compute the local dimension of
the Riesz product $P_{\theta}$ and we will apply Billingsley's
theorem to conclude Theorem \ref{Main}.

 For each $x\in
\mathbb{D}$ and $n\geq 1$, let
$$I_n(x)=I(x_1,\cdots,x_n)=\set{y\in \mathbb{D}:
y_k=x_k\textrm{ for }1\leq k\leq n}.$$ It is the $n$-cylinder
containing $x$, a ball of diameter $2^{-n}$. By the formula
(\ref{experence}), for any $n\geq \ell$, we have
    \begin{equation*}\label{measure-I_n}
P_{\theta}(I_n(x))=\frac{1}{2^n}\prod_{k=1}^{\lfloor
n/{\ell}\rfloor}\bigl(1+\theta \xi_k(x)\bigr).
    \end{equation*}
%By the Taylor expansion of $\log(1+\theta t) $ and the fact that
%$\xi_k^{2n}=1 $ and $\xi_k^{2n-1}=\xi_k$, we have
Recalling that $\xi_k(x)=+1$ or $-1$ for all $x$, by Taylor formula, we have
    \[  \log(1+\theta \xi_k(x)) = - \sum_{n=1}^{\infty}
  \frac{\theta^{2n}}{2n} +\sum_{n=1}^{\infty}
  \frac{\theta^{2n-1}}{2n-1} \xi_k(x).\]
Then for all points $x\in B_{\theta}$,
  \begin{equation*}
 \lim_{m\to\infty}\frac{1}{m}\sum_{k=1}^{m} \log(1+\theta \xi_k(x))
  =- \sum_{n=1}^{\infty} \frac{\theta^{2n}}{2n} + \sum_{n=1}^{\infty}
  \frac{\theta^{2n-1}}{2n-1}  \theta.
\end{equation*}
The right hand side can be written as $$
 \theta \log (1+\theta)-\frac{\theta-1}{2}\log
  (1-\theta^2)\\
  =\left[1-  H\left(\frac{1+\theta}{2}\right)\right]\log 2.$$
  It then follows that for
all points $x\in B_{\theta}$,
\begin{eqnarray*}
\lim_{n\to\infty}\frac{\log P_{\theta}(I_n(x))}{\log |I_n(x)|}
=\lim_{n\to\infty} \frac{\sum_{k=1}^{\lfloor n/{\ell}\rfloor}\log
\bigl(1+\theta \xi_k(x)\bigr)-\log 2^{n}}{\log 2^{-n}}=
1-\frac{1}{\ell}+\frac{1}{\ell}H\Bigl(\frac{1+\theta}{2}\Bigr).
\end{eqnarray*}
The proof is completed by applying Billingsley's theorem (\cite{Bi}). $\Box$

%\section{Proof of Theorem \ref{Main-HD-same}}
%Denote by $s$ the Hausdorff dimension of $A_0$ and by $t$ that of $X_0$. Suppose that $\{U_i\}$ is a covering of the set $A_0$
%such that $|U_i|\leq 2^{-n}$ for each $i$.

\medskip
\section{Proof of Theorem \ref{Main-box}}

It is clear that %the box dimension of $X_0$ is equal to
\[\dim_B X_0=\lim_{n\to\infty}\frac{\log_2N_n}{n}\]
if the limit exists, where $N_n$ is the cardinality of the following
set
\[
 \{(x_1x_2\cdots x_{n}): x_{\ell}x_{2{\ell}}=0 \ \text{ for } \ell\geq 1 \text{ such that }   \ 2{\ell}\leq
  n\}.
\]
Each equality $x_{\ell}x_{2\ell}=0$ defines a condition on the
sequence $(x_1\cdots x_n)$ which determines the cylinder
$I(x_1,\cdots, x_n)$. We observe that all these conditions can be
divided into ``independent" groups of conditions. Let
%
% To calculate $N_n$, we write the $n$ numbers $\{1,2, \dots,
%n\}$ in a form of a matrix, with each row having different number of
%$2$-divisors. That means we write the numbers in the form of row
%like
\begin{align*} C_0:&=\left\{1,\ 3,\ 5, \ \dots, \  2
n_0-1 \right\}, \\
C_1:&=\left\{2\cdot 1, \ 2\cdot 3, \ 2\cdot 5, \ \dots, \
2\cdot\big(2n_1-1\big) \right\}, \\
&\dots\\
C_k:&=\left\{2^k \cdot 1, \ 2^k\cdot 3 , \ 2^k\cdot 5, \ \dots, \
2^k\cdot (2n_k-1) \right\}, \\
&\dots\\
C_m:&=\left\{2^m \cdot 1\right\},
\end{align*}
where $n_k$ is the biggest integer such that $$2^k(2n_k-1)\leq n,
\quad \text{i.e.}, \quad
n_k=\left\lfloor\frac{n}{2^{k+1}}+\frac{1}{2}\right\rfloor$$ and $m$
is the biggest integer such that
$$2^m\leq n, \quad \text{i.e.}, \quad m=\lfloor\log_2 n\rfloor.$$
We have the decomposition $\{1, \cdots, n\}=C_0 \sqcup C_1\sqcup
\cdots\sqcup C_m$ and \[n_0>n_1>\cdots > n_{m-1}>n_m=1.\] The
conditions $x_{\ell}x_{2\ell}=0$ with $\ell$ in different columns in
the table defining $C_0, \cdots, C_m$ are independent. We are going
to use this independence to count the number of possible choices for
$(x_1,\cdots, x_n)$.

We have $n_m(=1)$ columns each of which has $m+1$ elements. Then we have
$a_{m+1}$ choices for $x_{\ell}$ with $\ell$ in the first column
since $(x_{\ell}, x_{2\ell})$ is conditioned to be different from  $(1,1)$. Each of the next
$n_{m-1}-n_{m}$ columns has $m$ elements, then we have $a_{m}^{n_{m-1}-n_m}$
choices for the $x_{\ell}$'s with $\ell$ in these columns. By
induction, we get
\[N_n=a_{m+1}^{n_m}a_{m}^{n_{m-1}-n_m}a_{m-1}^{n_{m-2}-n_{m-1}}\cdots a_1^{n_0-n_1}.\]

%
% By marking each number $\ell$ the value $x_{\ell}$, we know that the condition $x_{\ell}x_{2{\ell}}=0$ is equivalent to say in
%each column of the matrix, there is no consecutive ones. The number
%of possibilities for this case to happen is determined by the
%Fibonacci numbers. Let  $\set{a_k}$ be the Fibonacci sequence
%defined by
% $$a_0=1, \ a_1=2,\ a_2=3, \ a_3=5, \ a_4=8, \ \cdots, \ a_{k+2}=a_{k+1}+a_{k}.$$
%Precisely, the first $k$ numbers of each column has $a_{k}$
%possibilities to have no consecutive ones. Now we start to mark
%zeros and ones row by row. Consider the first row, we can freely
%mark zeros and ones, then we have
%$$2^{[\frac{n}{2}+\frac{1}{2}]}=\big(\frac{a_1}{a_0}\big)^{[\frac{n}{2}+\frac{1}{2}]}$$ possibilities.
%To pass from the first row to the second row, since we have to avoid
%the consecutive ones, we have to multiply the factor
%$$\big(\frac{a_2}{a_1}\big)^{[\frac{n}{2^2}+\frac{1}{2}]}.$$
%Go on this procedure, at the end, we obtain that
%\[
%  N_n= \prod_{k=0}^{[\log_2 n]}
%  \left(\frac{a_{k+1}}{a_k}\right)^{[\frac{n}{2^{k+1}}+\frac{1}{2}]}.
%\]
Now, the box dimension of the set $X_0$ equals to
\begin{align*}
  \lim_{n\to \infty}\frac{\log_2 N_n}{n}&= \lim_{n\to \infty} \frac{1}{n}\left(n_m\log_2
  a_{m+1}+\sum_{k=0}^{m}(n_{k-1}-n_k)\log_2a_k \right)\\
  &= \lim_{n\to \infty}\frac{1}{n}\left(\log_2
  a_{m+1}+\sum_{k=0}^{\lfloor\log_2 n\rfloor}\left(\left\lfloor\frac{n}{2^{k}}+\frac{1}{2}\right\rfloor-
  \left\lfloor\frac{n}{2^{k+1}}+\frac{1}{2}\right\rfloor\right)\log_2a_k
  \right)\\
    &=\sum_{k=1}^{\infty} \frac{\log_2 a_k}{2^{k+1}}.
\end{align*}

%Let $N_n$ be the cardinality of the following set
%\[
%  \{(x_1x_2\cdots x_{2n}): x_kx_{2k}=0 \ \text{ for} \ 1\leq k\leq
%  n\}.
%\]
%Then we have $N_{2n-1}=3N_{2n-2}$ and for $m=2^k(2n-1)$ with
%$k,n\geq 1$,
%\[N_m=\frac{2a_{k+1}}{a_k}N_{m-1},\]
%with $\set{a_k}$ being the Fibonacci sequence defined by
%  $$a_1=2,\ a_2=3, \ a_3=5, \ a_4=8, \ \cdots, \ a_{k+2}=a_{k+1}+a_{k}.$$
%The box dimension of $X_0$ is then determined by
%\[
%\lim_{n\to \infty} \frac{N_n}{2n\log 2}=\frac{1}{2\log
%2}\sum_{n=1}^{\infty}\frac{\log a_n}{2^n}=0.8242936\cdots.
%\]

\bigskip
\thanks{{\em Acknowledgement.}\quad  
%The authors wish to thank the Morning-side Center of
%Mathematics for the hospitality during their stay in the Summer of
%2009. 
This work is partially supported by NSFC10771164 (Ji-Hua Ma) and
NSFC10901124 (Lingmin Liao).}

\end{document}